\documentclass[a4paper,11pt]{article}
\pdfoutput=1
\usepackage[english]{babel}
\usepackage{fge,amssymb,amsthm,amsmath} 
\usepackage{mathrsfs} 
\usepackage{times} 
\usepackage{booktabs} 
\usepackage{multirow} 
\usepackage{subfig} 
\usepackage{xfrac} 
\usepackage{color}
\usepackage{authblk} 
\usepackage{a4wide} 

\newcommand{\Bd}{\mathcal{B}}
\newcommand{\Md}{\mathcal{M}}
\newcommand{\Id}{\mathrm{I}}
\newcommand{\lap}{\bigtriangleup}
\newcommand{\Qd}{\mathscr{Q}}
\newcommand{\R}{\mathbb{R}}
\newcommand{\X}{\mathrm{X}}
\newcommand{\Dsc}{\mathrm{D}}
\newcommand{\y}{\hat{y}}
\newcommand{\x}{\mathrm{x}}
\newcommand{\dat}{\mathrm{d}}
\newcommand{\Hv}{\mathrm{H}}
\newcommand{\Wv}{\mathrm{W}}
\newcommand{\La}{\mathcal{L}}
\newcommand{\setm}{\!\mathbin{\fgebackslash}\!}
\providecommand{\keywords}[1]{\textit{Keywords:} #1}

\makeatletter
\def\ifPositive#1{%
	\@ifnextchar{-}%
		{\expandafter\@secondoftwo\remove@to@nnil}%
		{\expandafter\@firstoftwo\remove@to@nnil}%
			#1\@nnil
}

\def\removesymb#1{\justlast#1\relax}
\def\justlast#1#2\relax{#2}
\makeatother

\newcommand{\NSc}[2]{ 
	 
	\ifPositive{#2} 
		{\ensuremath{#1\text{\sc{e}}#2}}
		{\ensuremath{#1\text{\sc{e-}}\removesymb{#2}}} 
}

\newcommand{\TimeFormat}[2]{ 
	\ensuremath{#1\text{m} #2\text{s}}
}

\begin{document}

\title{Radial basis function methods for optimal control of the convection-diffusion equation\footnote{Preprint}}
\author[1]{Pedro Gonz\'alez Casanova\footnote{Corresponding author, email address: \texttt{casanova@matem.unam.mx}}}
\author[1]{Jorge Zavaleta}
\affil[1]{\small Instituto de Matem\'aticas \\ Universidad Nacional Aut\'onoma de M\'exico \\ Cuidad Universitaria \\ Ciudad de M\'exico 04510 \\ MEXICO}
\date{}

\maketitle

\begin{abstract}
	PDE-constrained optimization problems have been barely solved by radial basis functions (RBFs) methods \cite{Pearson2013}. It is well known that RBF methods can attain an exponential rate of convergence when $C^{\infty}$ kernels are used, also, these techniques, which are truly scattered, are known to be flexible to discretize complex boundaries in several dimensions. On the other hand, exponential convergence implies an exponential growth of the condition number of the Gram matrix associated with these meshfree methods and global collocation techniques are known to be computationally expensive. In this paper, and in the context of optimal constrained optimization problems, we aim to explore a possible answer to both problems. Specifically, we introduce two local RBF methods: LAM-DQ based in the combination of an asymmetric local method (LAM), inspired in local Hermite interpolation (LHI), with the differential quadrature method (DQ), and LAM-LAM which use two times the local asymmetric method. The efficiency of these local methods against global collocation by solving several synthetic convection-diffusion control problems is analyzed. In this article, we also propose a preconditioning technique and treat the ill-conditioning problem by using extended arithmetic precision. We think that these local methods, which are highly parallelizable, shows a possible way to solve massive optimization control problems in an efficient way.
\par\vspace{0.5em}\noindent
\end{abstract}

\section{Introduction}

	Several works have appeared in the literature which deals with the analysis and formulation of numerical methods for the solution of distributed control problems in two or three dimension \cite{ZhouZJ2014}. These works have been formulated within the context of two general frames: the discretized-optimized and/or the optimize-discretized approaches. In particular, Galerkin methods have been proposed and analyzed within both frames, (see \cite{ZhouZJ2014} and references therein). On the other hand, we could only find one article in the literature, which uses radial basis function methods for constrained optimization problems, see \cite{Pearson2013}. More precisely in \cite{Pearson2013} the author solves  Poisson constrained optimization problems by using global RBFs symmetric and asymmetric collocation techniques.
	
	It is well known that a major limitation of global RBFs collocation techniques is that as the number of nodes or the shape parameter increases, the condition number of the corresponding Gram matrix grows. In the case of infinitely differentiable RBFs, convergence can be exponential but the corresponding condition number also increases in an exponential way \cite{Sarra2011}. 

	The current article is formulated within the context of the optimize-then-discretize approach and has the following objectives: 
\begin{itemize}
	\item Extend the application of RBFs meshfree methods to solve convection-diffusion constrained optimization problems.
	\item Use local RBFs methods, or more precisely, introduce a local asymmetric version, LAM, inspired on the local Hermite interpolation technique, LHI, \cite{Stevens2010}.
	\item Show that the proposed local methods can attain the same accuracy but with the advantage of a considerable reduction of the computing time, giving the possibility to extend its application for massive problems.
	\item Introduce a simple but effective preconditioner to invert the local matrices of the LAM method.
\end{itemize}

	We find that the local interpolation method, LHI, applied to the convection-diffusion control problem gives rise to a saddle point problem which is well known, in general, to be ill-posed. Although we manage to prove that under certain conditions the saddle point matrix is invertible, the result is not strong enough and we thus formulated a different approach. Specifically, we decouple the Euler Lagrange system of equations, corresponding to the control problem, by formulating a biharmonic problem which let us compute the state variable through LAM, and once this is done we use the state variable to calculate the control by DQ or using again LAM.
 
	It is worth to mention that, an important alternative to the ill-conditioned problem of RBFs collocation methods, is that they can be solved by using domain decomposition methods, see, \cite{Gonzalez2009}, \cite{Chinchapatnam2007}, \cite{ZhouX2003}. In this work, however, we are interested in comparing the results of the local methods with those obtained by solving the global asymmetric collocation (AC) method by using direct solvers with quad precision which is a current alternative that has been used by Kansa \cite{Kansa2017} and Sarra, \cite{Sarra2011} among others.

This paper is organized as follows. In section 2, we briefly state the continuous control problem and refer the reader to the proper references. Section 3, is devoted 
 to formulate the LAM-LAM and LAM-DQ local methods for the solution of the convection-diffusion constrained optimization problems as well as the preconditioning technique. Section \ref{NumExamples}, we present numerical examples to show the capabilities of the local methods. In section 5, conclusions are presented.
 
\section{The convection-diffusion control problem}

	Throughout this paper, we will be concerned with the solution of the following distributed control problem
\begin{equation}\label{Functional}
		\begin{array}{c}
			\min_{y,u} {\frac{1}{2}\Vert y-\y \Vert^2_{L^{2}(\Omega)} + \frac{\beta}{2}\Vert u 				\Vert^2_{L^{2}(\Omega)}}\\
			\\
			\mbox{s.t. } \mathcal{E} y=u\mbox{ in }\Omega,\quad\Bd y=g\mbox{ on }\partial\Omega
		\end{array}
	\end{equation}
where, $y$ is the state, $u$ the control, $\y$ the objective state, $\beta>0$ the penalty constant, ${\mathcal{E}}$ is a PDE stationary linear operator and $\Bd$ a Dirichlet, Neumann or Robin, boundary operator. Such problems were introduced by Lions in  \cite{Lions1968}.

	The distributed control problem (\ref{Functional}) can be equivalently formulated as a functional that incorporates the PDE constraints by means of Lagrange multipliers, (see \cite{Rees20101}), namely as
\begin{equation}\label{minimization_problem}
	\La(y,u,p_{1},p_{2}) = \frac{1}{2}\Vert y - \y \Vert^2_{L^{2}(\Omega)} + \frac{\beta}{2}\Vert u \Vert^2_{L^{2}(\Omega)} + \int_{\Omega}\left(\mathcal{E} y - u\right)p_{1} + \int_{\partial\Omega}\left(\Bd y - g\right)p_{2}.
\end{equation}

	Taking the Frechet derivative of the functional  (\ref{minimization_problem}) with respect to $y$, $u$ and $p$ it is possible to obtain the following Euler Lagrange equations in terms of the state $y$ and the control variable $u$,

\begin{equation}\label{PDE_System}
		\left.
		\begin{array}{cl}
			\mathcal{E} y = u & \mbox{ in }\Omega\\
			\Bd y = g       & \mbox{ on }\partial\Omega
		\end{array}\quad\right\vert\quad
		\begin{array}{cl}
			\beta\mathcal{E}^{*}u = \y-y & \mbox{ in }\Omega\\
			u = 0                      & \mbox{ on }\partial\Omega
		\end{array}
\end{equation}
where we have eliminated $p$ by using the equation $p = \beta u$. Here, $y$ and $u$ satisfying \eqref{PDE_System} are known as the optimal state and optimal control, respectively.

\section{Numerical schemes}

	The following schemes will be discretized by using multiquadric RBFs \emph{i.e.} $\Phi(x)=\sqrt{c+\Vert x\Vert^{2}}$, where $c$ is the shape parameter. We first describe the global asymmetric collocation and local methods to solve the minimization problem.
		
\subsection{Asymmetric collocation}
	
	In order to formulate the global asymmetric collocation scheme for the former coupled pair of equations \eqref{PDE_System} we first define the following ansatz
\[
		y(x)=\Hv(x)\lambda, \quad u(x)=\Hv(x)\mu,
\]
where $\Hv$ is known as the reconstruction vector, taken here as usual as
\[
		\Hv(x) = \left[\left.\underset{1\leq i\leq n} {\Phi\left(x-\x_{i}\right)}\;\right\vert\underset{1\leq\ell\leq n_{p}}{p_{\ell}(x)}\right]\in\R^{n+n_{p}},
\]
with $n$ the number of nodes and $n_{p}$ the number of polynomial terms. Taking the first $n_b < n$ nodes to be the boundary nodes, the resulting system of linear equations is given by
\[
	\begin{bmatrix}
		 G^{\Bd } & \beta E_{*}\\
		 -E       & G
	\end{bmatrix}
	\begin{bmatrix}
		\lambda \\ 
		\mu
	\end{bmatrix}=
	\begin{bmatrix}
		\dat \\ 
		0
	\end{bmatrix},
\]
where
\[
	E = \begin{bmatrix}
		0                        & 0                     \\
		\mathcal{E}\Phi_{\Omega} & \mathcal{E}P_{\Omega} \\
		0                        & 0
	\end{bmatrix},\quad
	E_{*} = \begin{bmatrix}
		0                            & 0                         \\
		\mathcal{E}^{*}\Phi_{\Omega} & \mathcal{E}^{*}P_{\Omega} \\
		0                            & 0
	\end{bmatrix},\quad
	G^{\Bd } = 
	\begin{bmatrix}
		\Bd \Phi_{\partial\Omega} & \Bd P_{\partial\Omega} \\
		\Phi_{\Omega}             & P_{\Omega}             \\
		P^{t}                     & 0
	\end{bmatrix},
\]
are square matrices of size $(n+n_{p})\times (n+n_{p})$, and $(\Bd \Phi_{\partial\Omega})_{j,i} = \Bd \Phi(x_{j}-x_{i})$, $(\Bd P_{\partial\Omega})_{j,\ell} = \Bd p_{\ell}(x_{j})$, $(\Qd \Phi_{\Omega})_{k,i} = \Qd\Phi(x_{k}-x_{i})$, $(\Qd P_{\Omega})_{k,i} = \Qd p_{\ell}(x_{k})$, for $\Qd = \mathcal{E}^{*},\mathcal{E},I$, with $I$ the identity operator, $G := G^{\Bd}$ is the standard Gram matrix for $\Bd  = I$, $P^t=\begin{bmatrix} P_{\partial\Omega}^{t} & P_{\Omega}^t \end{bmatrix}$ and
\[
	\dat = \left[\left.\underset{1\leq j\leq n_{b}} {g(\x_{j})}\;\right\vert\left. \underset{n_{b}+1\leq k\leq n} {\y(\x_{k})}\;\right\vert \underset{1\leq\ell\leq n_{p}}{0}\right]^{t}.
\]

	If $G^{\Bd } = G$, \emph{i.e.} taking $\Bd = I$, we solve this system through block LU factorization as follows:
\[
	\begin{bmatrix}
		G  & \beta E_{*} \\
		-E & G
	\end{bmatrix} =
	\begin{bmatrix}
		I_{n+n_{p}} & 0           \\
		-EG^{-1}    & I_{n+n_{p}}
	\end{bmatrix}
	\begin{bmatrix}
		G & \beta E_{*} \\
		0 & R
	\end{bmatrix},
\]
where $R=G+\beta EG^{-1}E_{*}$, the Schur complement of $G$ and $I_{n}$ is the identity matrix of size $n$.

\subsection{A local asymmetric scheme}

	The system of equations (\ref{PDE_System}) can be easily shown to be equivalent to the following boundary value elliptic problem, assuming $y$ is smooth enough
\begin{equation}\label{PDE_LAM}
	\begin{array}{rl}
		\Md y =\y         & \mbox{ in }\Omega         \\
		\mathcal{E} y = 0 & \mbox{ on }\partial\Omega \\			
		\Bd y = g         & \mbox{ on }\partial\Omega
	\end{array}
\end{equation}
where the differential operator $\Md$ is given by $\Md=\Id+\beta\mathcal{E}^{*}\mathcal{E}$. Although system (\ref{PDE_System}) can be directly discretized, it involves the solution of a saddle point problem, which is well known to be singular unless special conditions, e. g. inf-sup conditions, in the case of finite elements, are imposed. We thus find it more convenient to use system (\ref{PDE_LAM}) to compute the numerical solution.

	To formulate the LAM scheme of the system (\ref{PDE_LAM}) we consider the following: Let $\X\subset {\bar{\Omega}}$  be a set of $n$ scattered nodes and let $\X_{c}$ be a subset of $n_{c}$ nodes. Consider neighborhoods $\Dsc_{k}$ (\emph{e.g.} a disc of fixed radius) around the $k$-th point of $\X_{c}$ and label the nodes of $\Dsc_{k}\cap\X$ so that:
\begin{itemize}
	\item There are $n^{(k)}$ nodes in $\Dsc_{k}$ \emph{i.e.} $n^{(k)}=\#(\X\cap\Dsc_{k})$.
	\item The first node, $\x_{1}^{(k)}$ is the center of $\Dsc_{k}$.
	\item The first $n_{c}^{(k)}$ nodes are centers of other discs, \emph{i.e.} $n_{c}^{(k)}=\#(\X_{c}\cap\Dsc_{k})$.
	\item The following $n_{b}^{(k)}$ nodes lie on $\partial\Omega$, $n_{b}^{(k)}=\#(\partial\Omega \cap(\Dsc_{k}\setm \X_{c}))$.
	\item The remaining $n_{\iota}^{(k)}$ nodes belong to the interior of $\Omega$ (and none of them are centers of any disc), so that $n^{(k)}=n_{c}^{(k)}+n_{b}^{(k)}+n_{\iota}^{(k)}$.
\end{itemize}
	For each disk, the method forms a local system whose solutions are used to build a global sparse matrix. The solution of this global system gives the approximated values of the PDE system \eqref{PDE_LAM} at the centers $\X_{c}\subset\Omega$.
	
	Choosing a conditionally positive definite radial basis function $\Phi$ of order $m$ and let $n_{p}$ be dimension of the corresponding polynomial space, we define the reconstruction vector
\begin{equation*}
	\Hv^{(k)}(x) = \left[\left.\underset{1\leq j\leq n^{(k)}} {\Phi\left(x-\x_{j}^{(k)}\right)}\;\right\vert\underset{1\leq\ell\leq n_{p}}{p_{\ell}(x)}\right]\in\R^{n^{(k)}+n_{p}}.
\end{equation*}

	Defining the following ansatz
 \[
 	y^{(k)}(x)=\Hv^{(k)}(x)\lambda^{(k)},
 \]
we obtain the local linear system
\begin{equation*}
	A^{(k)}\lambda^{(k)} =
	\begin{bmatrix}
		\Phi            & P             \\
		\Bd\Phi         & \Bd P         \\
		\mathcal{E}\Phi & \mathcal{E} P \\
		\Md\Phi         & \Md P         \\
		P^{t}           & 0
	\end{bmatrix}
	\lambda^{(k)}
	= \dat^{(k)}
\end{equation*}
with the data vector
{\small
	\begin{equation*}
		\left.\dat^{(k)} = \left[\underset{1\leq j\leq n_{c}^{(k)}} {y\left(\x_{j}^{(k)}\right)}\;\right\vert\left.
		\underset{n_{c}^{(k)}< j\leq n_{c}^{(k)}+n_{b1}^{(k)}}{g\left(\x_{j}^{(k)}\right)}\;\right\vert\left.
		\underset{n_{c}^{(k)}+n_{b1}^{(k)}< j\leq n_{c}^{(k)}+n_{b}^{(k)}}{0}\right.\left\vert\;
		\underset{n_{c}^{(k)}+n_{b}^{(k)}< j\leq n^{(k)}}{\y\left(\x_{j}^{(k)}\right)}\;\right\vert
		\underset{1\leq\ell\leq n_{p}}{0}\right]^{t}
	\end{equation*}}
where $n_{b1}^{(k)}$ and $n_{b2}^{(k)}$ are the number of boundary points for each of the boundary conditions, so that $n_{b}^{(k)}=n_{b1}^{(k)}+n_{b2}^{(k)}$. Solving for $\lambda^{(k)}$ we obtain the local solution 
\begin{equation}\label{weights}
	y^{(k)}(x)=\Hv^{(k)}(x)\left(A^{(k)}\right)^{-1}\dat^{(k)}=W^k(x)\dat^{(k)},
\end{equation}
where $\Wv^{(k)}$ is known as the vector of weights. Using this last expression it is possible to compute $\Qd u^{(k)}$ for any differential operator $\Qd$ through $\Qd u^{(k)}(x)=\left(\Qd \Wv^{(k)}\right)(x)\dat^{(k)}$.

	Denote by $y_{c}=\left[y\left(\x_{1}^{(k)}\right)\right]_{k=1}^{n_{c}}\in\R^{n_{c}}$ the vector of the values of $y$ at each of the centers. Then for each  $k$, the unknown elements of  $\dat^{(k)}$ belong to $y_{c}$.

Consider now the following system of equations
\begin{align}\label{weights2}
	\y\left(\x_{1}^{(k)}\right)=\Md y\left(\x_{1}^{(k)}\right)= \Md \Hv^{(k)}(\x_{1}^{(k)})\left(A^{(k)}\right)^{-1}\dat^{(k)} = \Wv_{\Md}^{(k)}\left(\x_{1}^{(k)}\right)\dat^{(k)}
\end{align} 
for $k=1,\ldots,n_{c}$ and $\Wv_{\Md}^{(k)}=\Md\Wv^{(k)}$. 

	This is a linear system in $y_{c}$, whose elements are the approximated solution of the PDE system (\ref{PDE_LAM}), at the centers, and which can be written as  $Sy_{c}=b$. Note that since in each $\dat^{(k)}$ there are only a few number of centers, \emph{i.e.} $n_{c}^{(k)}$ is relatively small, the matrix $S$ is sparse and thus standard preconditioning techniques can be used.

	In order to build the matrix $S$, we compute the weights by solving the following equation, (see equation (\ref{weights2})),
\begin{align}\label{weights3}
	\Wv_{\Md}^{(k)}\left(\x_{1}^{(k)}\right) = \Md \Hv^{(k)}(\x_{1}^{(k)})\left(A^{(k)}\right)^{-1}  
\end{align}

	Once the state $y$ has been computed, we can obtain the control $u$, through one of the following two algorithms:
\begin{enumerate}
	\item {\bf Local asymmetric method (LAM)}. Solve the problem for $u$ by means of,
	\begin{equation}\label{LAM-LAM}
		\begin{array}{cl}
			\beta\mathcal{E}^{*}u = \y-y & \mbox{ in }\Omega         \\
			u=0                          & \mbox{ on }\partial\Omega
		\end{array}
	\end{equation}
	using the computed values of  $y$.
	\item {\bf Differential quadrature (DQ)}. Where we evaluate  
	\[
		u = \mathcal{E} y
	\]
	by discretizing the operator using the differential quadrature technique.
\end{enumerate}

	We shall denote the first scheme by LAM-LAM and by LAM-DQ to the second one. We omit the description of the LAM-LAM algorithm, as the second part, the system (\ref{LAM-LAM}), has been essential already described. We thus briefly recall the differential quadrature method for this problem. 

	The  main point of the RBF differential quadrature method, see Shu \cite{Shu2003}, is to build a discrete operator $\tilde{\mathbf{E}}$ which approximates the continuous linear differential operator $\mathcal{E}$. Its construction can be summarized as follows. First, we solve the following system
\begin{equation} \label{DQ1}
	\mathcal{E}\Phi(x) \Big{|}_{x=\x_{k}}= \sum^{n_{k}}_{j=1}{ w^{\mathcal{E}}_{k,j}  \Phi(\x_{k,j}) }, \quad k=1,2,\ldots,N
\end{equation}
where the nodes $\{\x_{k,j}\}_{j=1}^{n_k}\subset\Omega$ are the $n_{k}$ nearest points to $\x_{k}\in \X\subset\Omega$. For simplicity, we have taken $\Phi$ to be a strictly positive definite radial basis function, (the formulation also holds for conditional positive radial basis functions). It is well known, see  \cite{Wendland2004}, that the system (\ref{DQ1}) is invertible. Once the coefficients $w^{\mathcal{E}}_{k,j}$ are computed, the approximated discretization of the operator $\mathcal{E}$ of a smooth enough function $u$ is given by
\[
	\mathcal{E} u(x)  \Big{|}_{x=\x_{k}}     \approx  \tilde{\mathbf{E}}   u(\x_{k}) = \sum^{n_k}_{j=1} { w^{\mathcal{E}}_{k,j} u(\x_{k,j}) }, \quad k=1,2,\ldots,N . 
\]

	Note that unlike LAM approach, see equation (\ref{weights3}), the differential quadrature technique does not include the boundary operator $\Bd$ in the computation of the weights, equation (\ref{DQ1}), see \cite{Shu2003}.

\section{Numerical examples}\label{NumExamples}
	
	In this section, we will discuss different examples to illustrate our main contribution. Specifically, that the proposed local algorithms can attain errors which are comparable to the global asymmetric colocation technique but a much lower computational cost. The analysis of the numerical experiments for these techniques is not trivial due to the existence of three parameters that simultaneously controls the quality of the results. These parameters are the fill distance, the penalty constant, $\beta$, and the shape parameter, $c$.

	The experiments were set up in the following way: given a total number of nodes $n$, we vary the values of $\beta$ and/or $c$, showing that we obtain completely different results. In fact, although the error can be good the condition number can be close to the machine precision, which means that we have problems that are numerically ill-posed and the result may not be reliable. On the other hand, we can have a good condition number, which means that the scheme is stable, but the error can be very poor. The goal then is to find the appropriate parameters that guarantee both stability and good numerical errors. To do this, we look for values of $\beta$ and $c$ for which the error is minimal and the condition number of the Gram matrix which is within the used precision. The reason for this criterion is that the condition number tells us, approximately, how many digits of the error are reliable. In other words, for a condition number of $10^k$, up to $k$ digits of accuracy may be lost within the floating-point arithmetic used. We remark that the computed condition number is only an approximation that serves as a bound for the exact value of the maximum inaccuracy that may occur in the algorithms. In the case of the LAM-LAM and LAM-DQ methods, the restriction of the condition number is imposed on the local Gram matrices.

	The numerical results obtained by each method are compared with respect to the value of $\beta$, although the values of $c$ do not necessarily have to coincide for both techniques since one method is global and the other local. The results are presented using multiquadric RBFs and quadruple precision to further investigate the performance of the methods as well as the effect of the condition number. Finally, independently of the values for $\beta$ and $c$, there are problems that requires a greater number of local nodes to obtain good numerical errors.
	
\subsection{Problem 1}

	The first test problem that we would like to analyze is Poisson control given by:
\[
	\begin{array}{rcl}
		-\lap y =u, & -\beta\lap u= \y - y              & \mbox{ in }\Omega            \\
		y=g,        & u=0                               & \mbox{ on }\partial\Omega    \\
		            &                                   &                              \\
	                & \y=\sin{\pi x_{1}}\sin{\pi x_{2}} &                              \\
	                & g=0                               &
	\end{array}
\]
with exact solution given by
\begin{align*}
	y_{\beta}(x_{1},x_{2}) &= \frac{1}{1+4\beta(\pi)^{4}}\sin{\pi x_{1}}\sin{\pi x_{2}}\\
	u_{\beta}(x_{1},x_{2}) &= \frac{2\pi^{2}}{1+4\beta(\pi)^{4}}\sin{\pi x_{1}}\sin{\pi x_{2}}.
\end{align*} 

	Since we want to restrict the values of the condition numbers $\kappa$, corresponding to the local scheme, we shall use the value of $\kappa=\max_{k} \kappa(A^{(k)})$ to measure the numerical ill-posedness, meanwhile for the global method we use $k=\kappa(G)$.

	Table \ref{tabla:ResP1} contains the values $\Vert y-\y\Vert_{L_{2}(\Omega)}$ for th state $y$ and $\Vert u\Vert_{L_{2}(\Omega)}$ for control $u$; the relative error, $RE_{y} = \Vert  y-y_{\beta}\Vert_{L_{2}(\Omega)}/\Vert  y_{\beta}\Vert_{L_{2}(\Omega)}$,   $RE_{u}= \Vert  u-u_{\beta}\Vert_{L_{2}(\Omega)}/\Vert  u_{\beta}\Vert_{L_{2}(\Omega)}$ respectively and  the Cost $=\sfrac{(\|y-\hat{y}\|^2_{L_{2}(\Omega)} + \beta\|u\|^2_{L_{2}(\Omega)})}{2}$, where $\Vert  f\Vert_{L_{2}(\Omega)}^{2}=\sum\limits_{k=1}^{n} \vert f(\x_{k})\vert^{2}$.
{\scriptsize
	\begin{table}[t!]
		\begin{center}
			\begin{tabular}{c c c c c c c}\toprule[1.3pt]
				            & \multicolumn{3}{c}{\bf LAM-DQ}                      & \multicolumn{3}{c}{\bf AC}\\ \cmidrule(lr){2-4} \cmidrule(lr){5-7}
				$\beta$     & $10^{-4*}$      & $10^{-6}$       & $10^{-10}$      & $10^{-4}$        & $10^{-6}$       & $10^{-10}$      \\
				$c$         & \NSc{6.00}{-01} & \NSc{1.00}{00}  & \NSc{1.00}{00}  & \NSc{3.00}{-01}  & \NSc{4.00}{-01} & \NSc{4.00}{-01} \\
				$RE_y$      & \NSc{4.30}{-06} & \NSc{2.94}{-07} & \NSc{2.59}{-11} & \NSc{7.15}{-09}  & \NSc{3.23}{-09} & \NSc{4.28}{-12} \\
				$RE_u$      & \NSc{3.04}{-04} & \NSc{8.65}{-05} & \NSc{8.31}{-07} & \NSc{6.33}{-09}  & \NSc{1.98}{-07} & \NSc{6.59}{-08} \\
				$y$         & \NSc{3.32}{-01} & \NSc{3.45}{-03} & \NSc{3.45}{-07} & \NSc{4.22}{-01}  & \NSc{4.39}{-03} & \NSc{4.39}{-07} \\
				$u$         & \NSc{1.68}{02}  & \NSc{1.75}{02}  & \NSc{1.75}{02}  & \NSc{2.14}{02}   & \NSc{2.22}{02}  & \NSc{2.23}{02}  \\
			    Cost        & \NSc{1.47}{00}  & \NSc{1.53}{-02} & \NSc{1.53}{-06} & \NSc{2.38}{00}   & \NSc{2.47}{-02} & \NSc{2.48}{-06} \\
			    $\kappa$    & \NSc{4.87}{26}  & \NSc{1.51}{24}  & \NSc{1.42}{24}  & \NSc{2.03}{24}   & \NSc{9.05}{26}  & \NSc{9.05}{26}  \\
				$\kappa(S)$ & \NSc{4.68}{07}  & \NSc{3.81}{04}  & \NSc{1.39}{00}  &                  &                 &                 \\ 
				Time        & \TimeFormat{00}{46} & \TimeFormat{00}{09} & \TimeFormat{00}{09} & \TimeFormat{01}{51} & \TimeFormat{01}{51} & \TimeFormat{01}{51} \\ \bottomrule[1.3pt]
			\end{tabular} \caption{Results from problem 1. For LAM-DQ $n^{(k)} = 50$, except for * where $n^{(k)} = 100$. In both cases $n = 622$.}\label{tabla:ResP1}
		\end{center}
	\end{table}}

\begin{figure}[t!]
   \begin{center}
   	  \subfloat{LAM-DQ, $n^{(k)} = 50$}\hspace{0.35\textwidth}
      \subfloat{AC\hspace{1cm}}\\
      \addtocounter{subfigure}{-2}
      \subfloat[$RE_y$]{\includegraphics[width=0.46\textwidth]{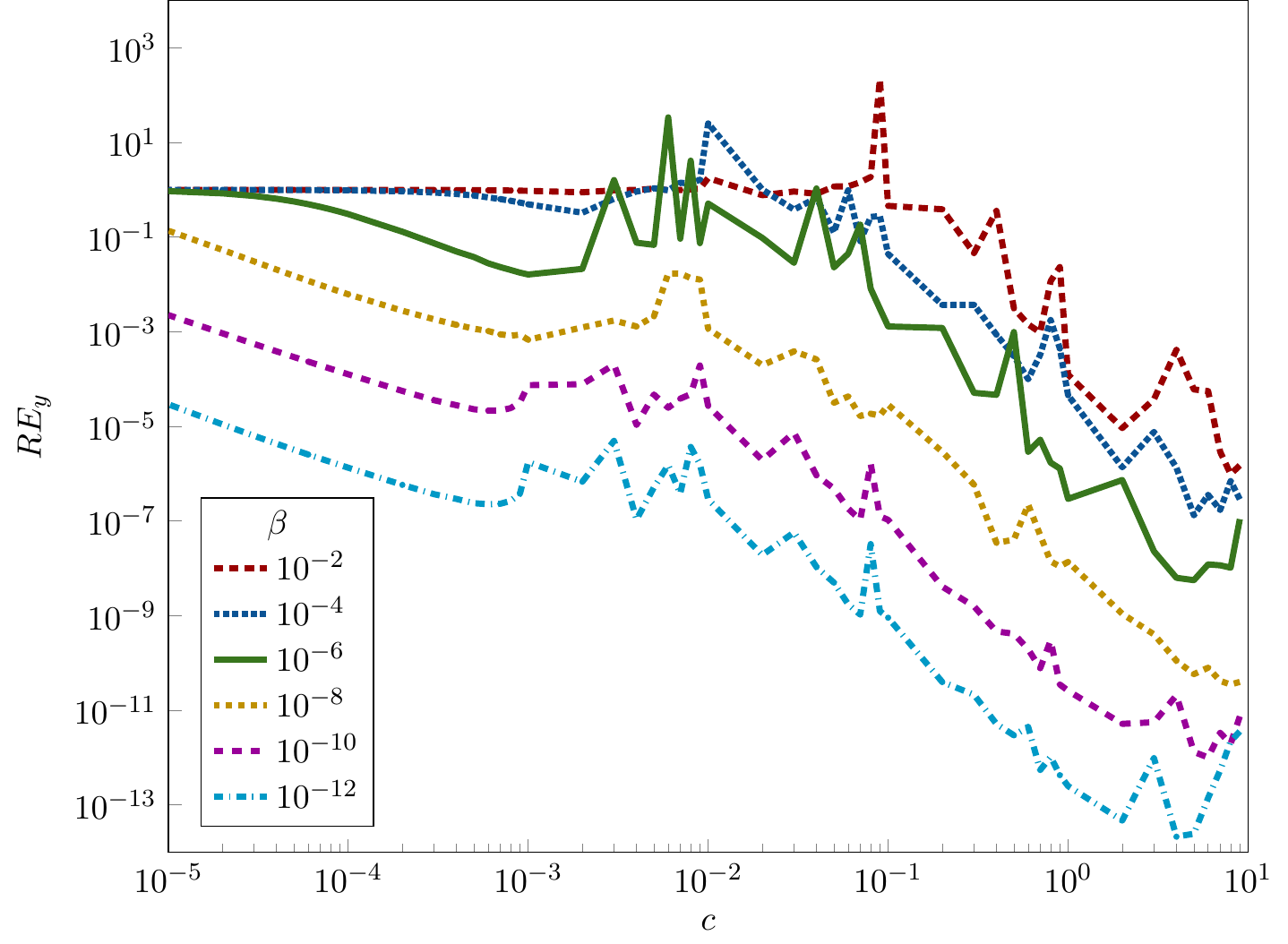}}\hspace{0.05\textwidth}
      \subfloat[$RE_y$]{\includegraphics[width=0.46\textwidth]{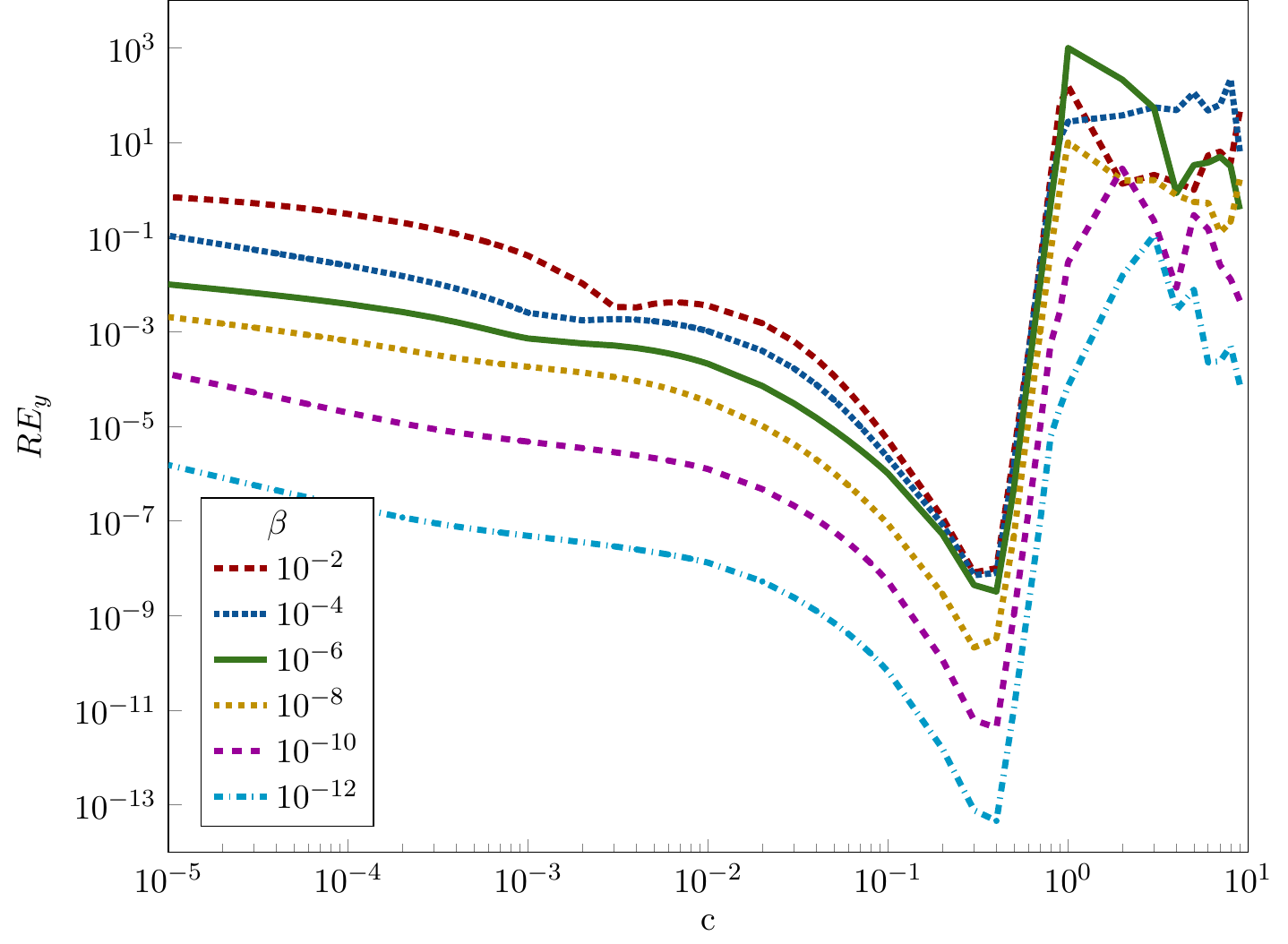}}\\
      \subfloat[$\kappa = \max_{k}\kappa(A^{(k)})$]{\includegraphics[width=0.46\textwidth]{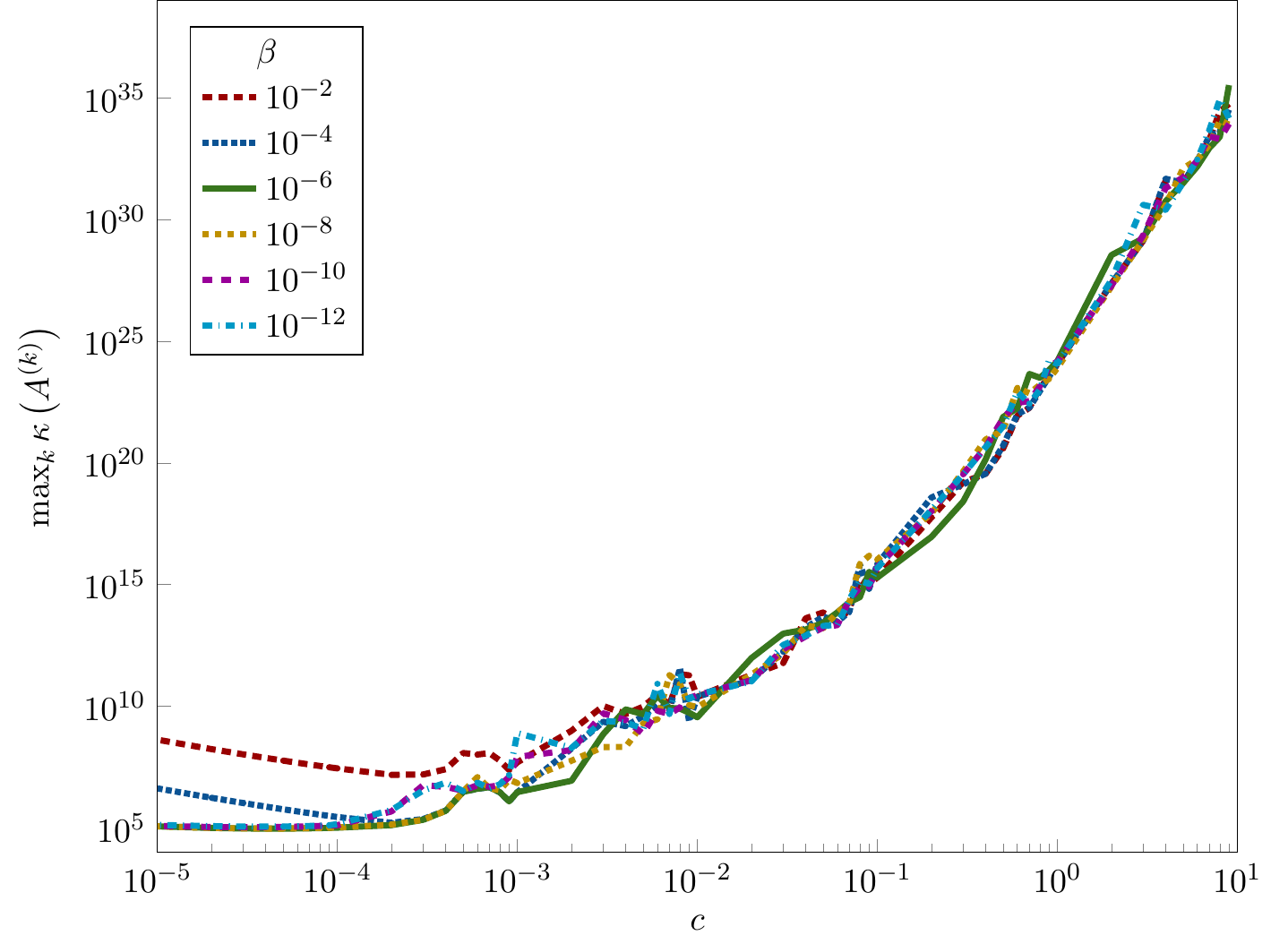}}\hspace{0.05\textwidth}
      \subfloat[$\kappa = \kappa(G)$, ]{\includegraphics[width=0.46\textwidth]{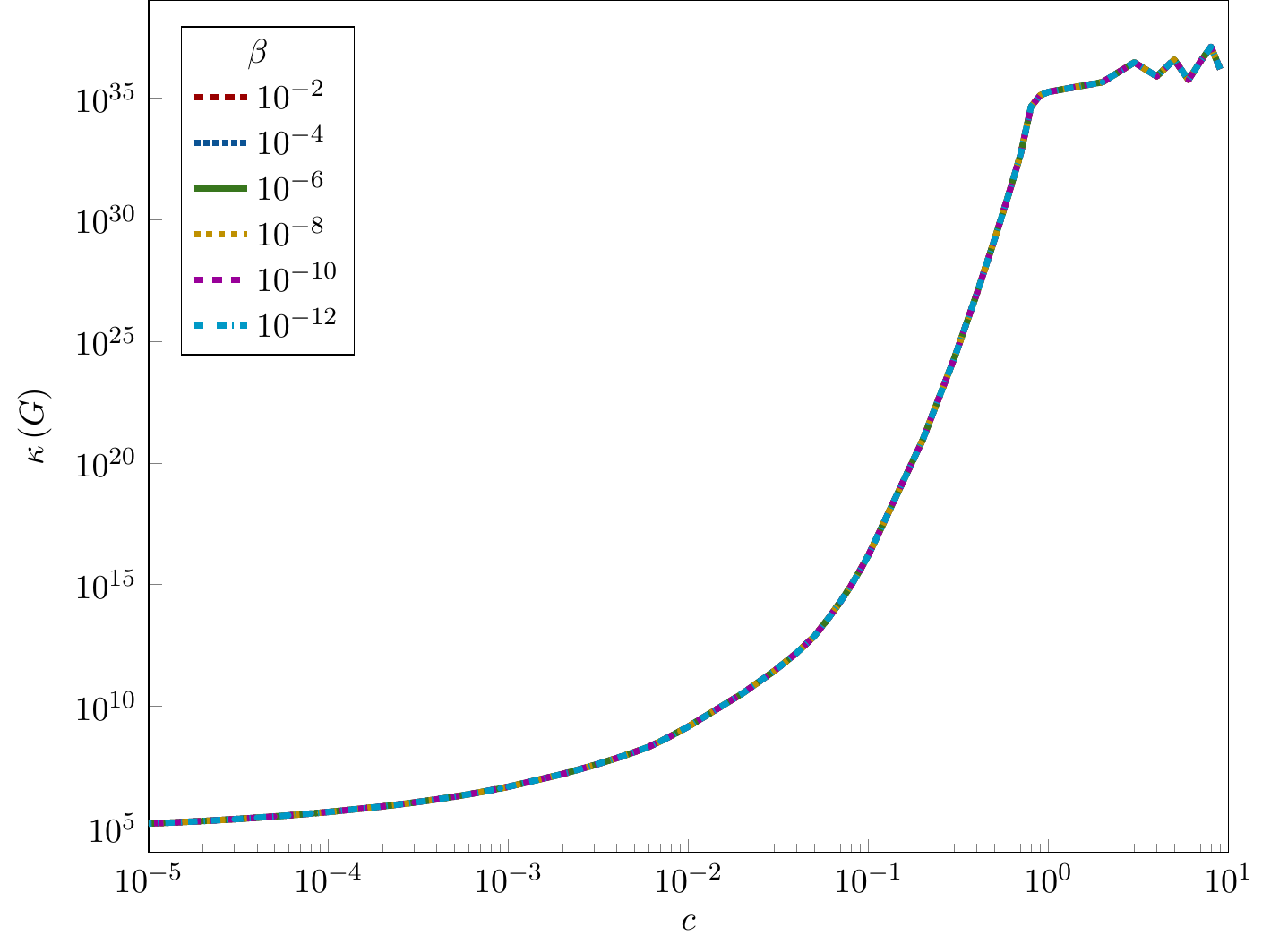}}
   \end{center}
   \caption{Comparison between the values of the relative error ($RE_{y}$) and the condition number ($\kappa$), by varying the shape parameter $c$. These calculations were obtained using quadruple precision, and different values of the penalty constant $\beta$.}\label{figura:CBetaP1QP}
\end{figure}

	From  table \ref{tabla:ResP1} we can observe that for small values of $\beta$ the errors obtained by global collocation and LAM-DQ techniques are comparable. Moreover, for large values of $\beta$, it is possible to change the number of nodes in the local systems to improve the LAM-DQ error. It is important to note that for small values of $\beta$ it is possible to obtain similar errors for both methods, but using a relatively small number of nodes for the local systems, which has a considerable impact on the computing time. Even when more nodes are used in local systems for large values of $\beta$, the computing time is still lower than the one used for AC. In addition, as $\beta\to0$, we have $\kappa(S)\to 1$, suggesting that the method is highly stable for these cases.

	Figure \ref{figura:CBetaP1QP} shows in detail the effect of the variation of $\beta$ and $c$ on the error and the condition number. We can see that as $\beta$ tends to zero and the value of $c$ increases the error decreases, so for both methods the results can be improved with respect to the error but they may be unreliable because of the conditioning when $c$ is increased. The results that are reported in the table \ref{tabla:ResP1} are far from the values of $\kappa$ for which the solutions are affected by rounding errors with respect to the precision used, still we obtain errors below $10^{-5}$. It is important to mention that as in our case, in \cite{Rees20102} the authors observe that as the value of $\beta$ decreases so does $\Vert y-\y\Vert_{L_{2}(\Omega)}$. In their case, it was only possible to explore this for values up to $\beta = 10^{-6}$, due to the limitation of their iterative methods designed for the finite element method.

\begin{figure}[t!]
   \begin{center}
      \subfloat[Relative error]{\includegraphics[width=0.46\textwidth]{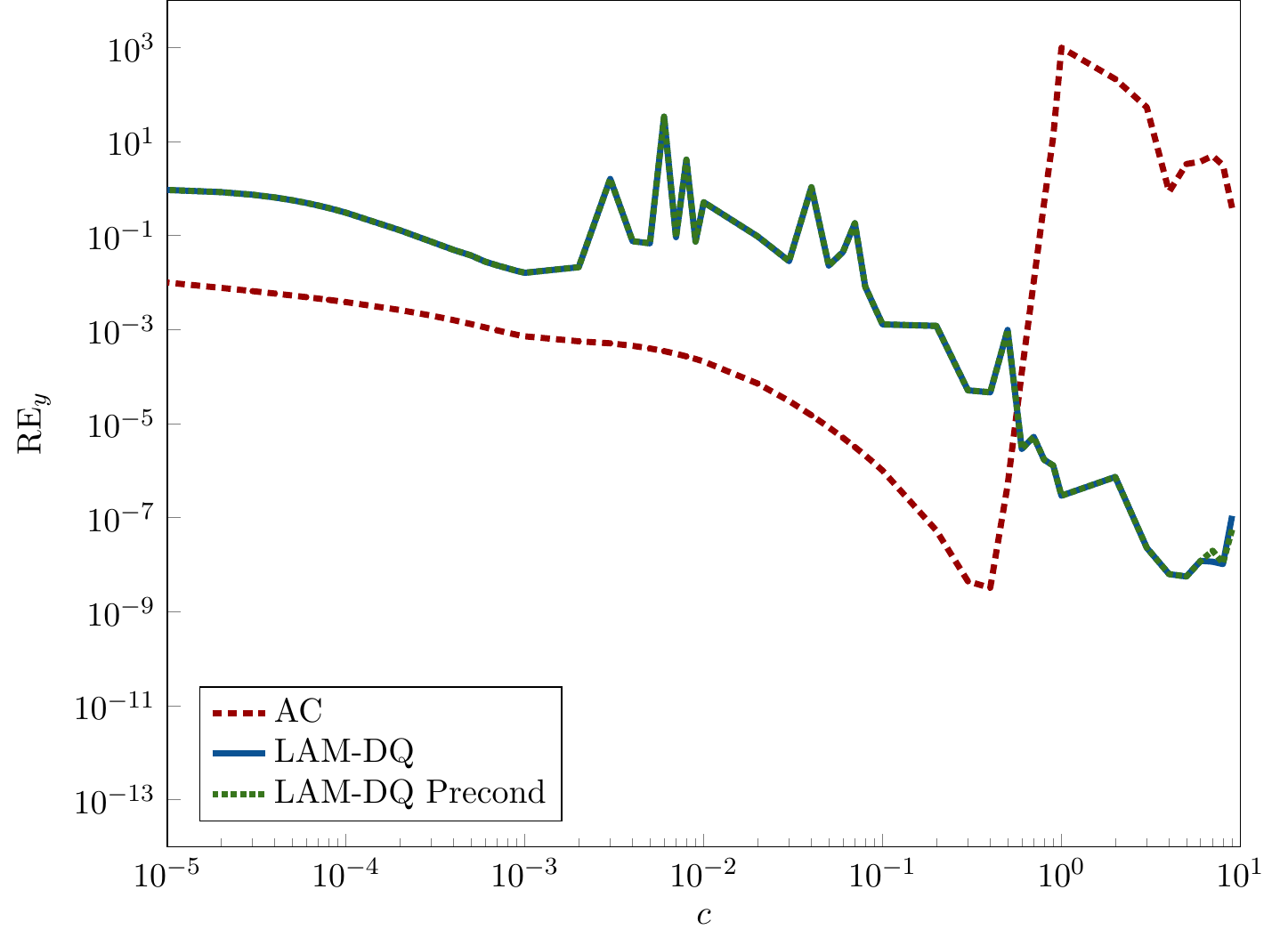}}\hspace{0.05\textwidth}
      \subfloat[Condition number]{\includegraphics[width=0.46\textwidth]{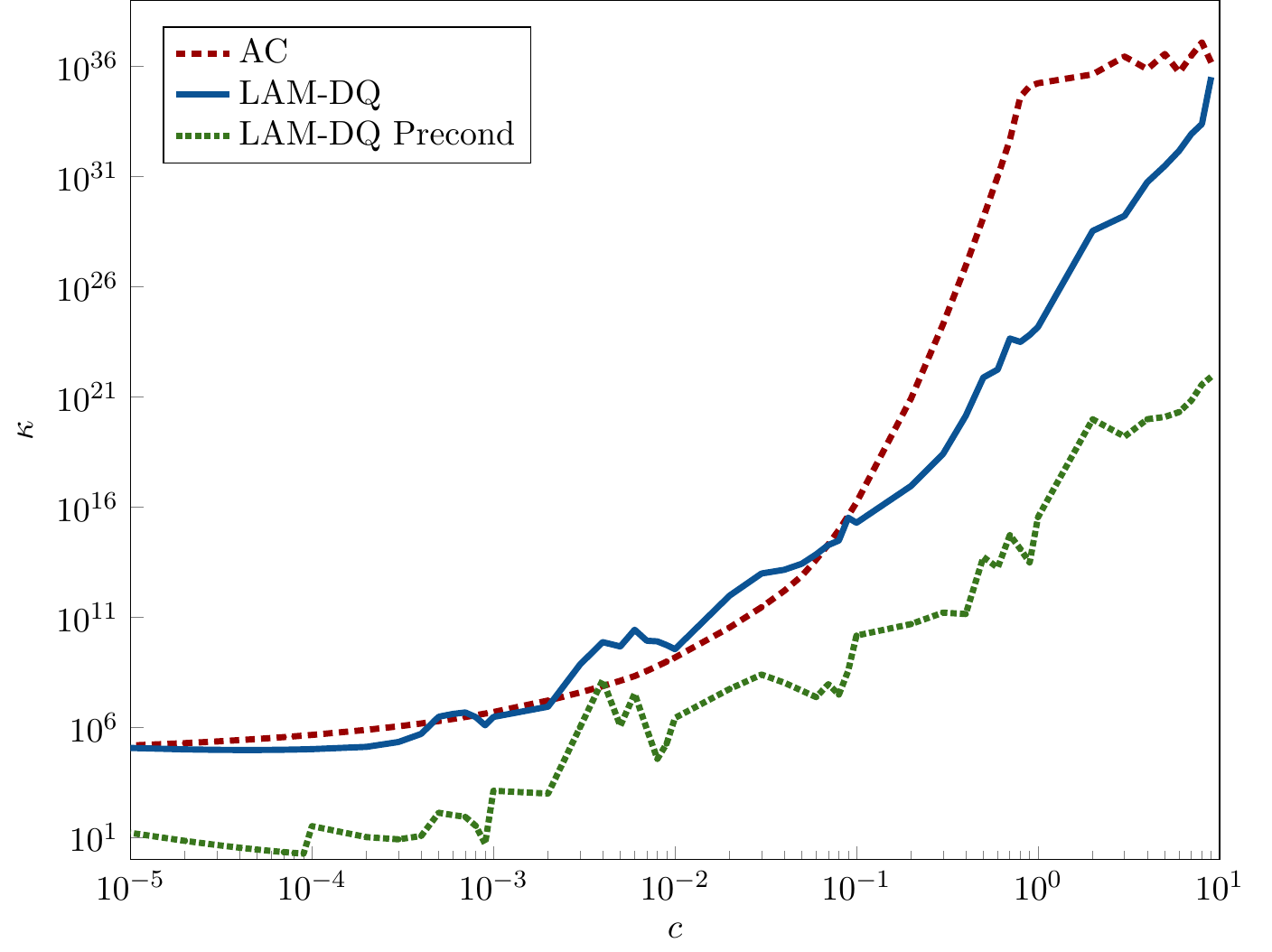}}
   \end{center}   
   \caption{Comparison of different methods for $\beta=10^{-6}$}\label{figura:CMetodoP1QP}
\end{figure}

	We also analyzed the use of a preconditioner for LAM-DQ, figure \ref{figura:CMetodoP1QP} shows a comparison of the methods for $\beta = 10^{-6}$. The point we want to emphasize here is that it is possible to reduce the conditioning of the local matrices $A^{(k)}$ in such a way that the results obtained for large values of $c$ are reliable. In this particular example, when using the preconditioner $P^{(k)}A^{(k)}$, with $P^{(k)} = (A^{(k)}_{*})^{-1}$, where $A^{(k)}_{*}$ is obtained in the same way as $A^{(k)}$ just by using the shape parameter $\hat{c}\neq c$, where $\hat{c} = c + \delta$ with $\delta$ small. For example, for the case of the figure \ref{figura:CMetodoP1QP}, given $c = m\times10^{\alpha}$ it was taken $\delta = 0.001\times10^{\alpha}$, such so that $\hat{c} = (m + 0.001)\times10^{\alpha}$, obtaining an error of the same size as in the case of LAM-DQ, but with a lower condition number, even reaching a difference up to 14 orders of magnitude for $c = 9$ where the condition number is around $10^{35}$ and in the case of AC up to 16 orders of magnitude for $ c = 8$ where the condition number reaches $10^{37} $, while the condition number for LAM-DQ Precond is around $10^{21} $ for both values of $c$.
	
{\scriptsize
	\begin{table}[t]
		\begin{center}
			\begin{tabular}{c c c c }\toprule[1.3pt]
				            & {\bf LAM-DQ}    & {\bf LAM-DQ Precond} & {\bf AC}\\ \cmidrule(lr){2-4}
				$\beta$     & $10^{-6}$       & $10^{-6}$       & $10^{-6}$       \\
				$c$         & \NSc{1.00}{00}  & \NSc{5.00}{00}  & \NSc{4.00}{-01} \\
				$RE_y$      & \NSc{2.94}{-07} & \NSc{5.63}{-09} & \NSc{3.23}{-09} \\
				$RE_u$      & \NSc{8.65}{-05} & \NSc{6.09}{-06} & \NSc{1.98}{-07} \\
				$y$         & \NSc{3.45}{-03} & \NSc{3.45}{-03} & \NSc{4.39}{-03} \\
				$u$         & \NSc{1.75}{02}  & \NSc{1.75}{02}  & \NSc{2.22}{02}  \\
			    Cost        & \NSc{1.53}{-02} & \NSc{1.53}{-02} & \NSc{2.47}{-02} \\
			    $\kappa$    & \NSc{1.51}{24}  & \NSc{1.25}{20}  & \NSc{9.05}{26}  \\
				$\kappa(S)$ & \NSc{3.81}{04}  & \NSc{1.91}{06}  &                 \\ 
				Time        & \TimeFormat{00}{09} & \TimeFormat{00}{20} & \TimeFormat{01}{51}  \\ \bottomrule[1.3pt]
			\end{tabular} \caption{Comparison of methods for $\beta = 10^{-6}$. For LAM-DQ and LAM-DQ Precond, $n^{(k)} = 50$. In all cases $n = 622$.}\label{tabla:CompMetodosP1QP}
		\end{center}
	\end{table}}	
	
	Table \ref{tabla:CompMetodosP1QP} compares the performance of LAM-DQ Precond with the best values reported in table \ref{tabla:ResP1} for $\beta = 10^{-6}$. In particular, for $c = 5$ and $\hat{c} = 5.001$, values obtained for $RE_y$ and $RE_u$ are almost of the same magnitude as for the global method, but with a lower condition number and still preserving a much lower computation time. There is clearly more room to improve this part, especially in the process of finding the optimal value of $\delta$ and thus looking for more efficient preconditioners.

\subsection{Problem 2}

	The next test problem is a convection-diffusion control problem for which there is no exact solution, given by
\begin{align*}		
	\begin{array}{rrl}
		(-\epsilon\lap+\omega\cdot\nabla)y=u, & \beta(-\epsilon\lap-\omega\cdot\nabla)u= y-\y & \mbox{ in }\Omega         \\
		y=g,                                  & u=0                                           & \mbox{ on }\partial\Omega \\
		                                      &
	\end{array}\\
	\begin{array}{rll}
		& \y=0                                                                                &\\
		& g = \begin{cases}
					1 & \mbox{in }\{0\}\times\left[\tfrac{1}{2},1\right]\cup[0,1]\times\{1\} \\ 
					0 & \mbox{elsewhere}
			  \end{cases}                                                                     & \\
		&  \omega=(\cos{\theta},\sin{\theta}),\mbox{ with } \theta=\displaystyle\frac{\pi}{6} & \\
		& \epsilon=\displaystyle\frac{1}{200}                                                 &
	\end{array}
\end{align*}

	This example corresponds to a boundary layer problem, which is of interest due to the sharp gradient attained at the boundary layer. Table \ref{tabla:ResP6} contains the values $\Vert y-\y\Vert_{L_{2}(\Omega)}$ for the state $y$ and $\Vert u\Vert_{L_{2}(\Omega)}$ for the control $u$.
	{\scriptsize
	\begin{table}[t!]
		\begin{center}
			\begin{tabular}{c c c c c c c}\toprule[1.3pt]
				            & \multicolumn{3}{c}{\bf LAM-DQ}                      & \multicolumn{3}{c}{\bf AC}\\ \cmidrule(lr){2-4} \cmidrule(lr){5-7}
				$\beta$     & $10^{-2}$       & $10^{-6}$       & $10^{-10}$      & $10^{-2}$        & $10^{-6}$       & $10^{-10}$      \\
				$c$         & \NSc{9.00}{-04} & \NSc{7.00}{-03} & \NSc{7.00}{-03} & \NSc{1.00}{-04}  & \NSc{4.00}{-04} & \NSc{4.00}{-04} \\
				$y$         & \NSc{3.72}{00}  & \NSc{1.31}{-02} & \NSc{1.32}{-06} & \NSc{3.94}{00}   & \NSc{1.10}{-02} & \NSc{1.10}{-06} \\
				$u$         & \NSc{2.64}{01}  & \NSc{2.96}{-01} & \NSc{2.97}{-05} & \NSc{3.84}{01}   & \NSc{2.61}{02}  & \NSc{2.61}{02}  \\
			    Cost        & \NSc{1.04}{01}  & \NSc{8.61}{-05} & \NSc{8.69}{-13} & \NSc{1.51}{01}   & \NSc{3.42}{-02} & \NSc{3.42}{-06} \\
			    $\kappa$    & \NSc{1.88}{05}  & \NSc{5.83}{08}  & \NSc{3.25}{08}  & \NSc{4.52}{05}   & \NSc{1.49}{06}  & \NSc{1.49}{06}  \\
				$\kappa(S)$ & \NSc{2.52}{04}  & \NSc{1.14}{00}  & \NSc{1.00}{00}  &                  &                 &                 \\ 
				Time        & \TimeFormat{00}{10} & \TimeFormat{00}{10} & \TimeFormat{00}{10} & \TimeFormat{01}{52} & \TimeFormat{01}{52} & \TimeFormat{01}{52} \\ \bottomrule[1.3pt]
			\end{tabular} \caption{Results from problem 2. For LAM-DQ $n^{(k)} = 50$. In both cases $n = 622$.}\label{tabla:ResP6}
		\end{center}
	\end{table}}	
	
\begin{figure}[t!]
   \begin{center}
   	  \subfloat{LAM-DQ, $n^{(k)} = 50$}\\
      \addtocounter{subfigure}{-1}
      \subfloat[State, $c = 10^{-8}$]{\includegraphics[width=0.45\textwidth]{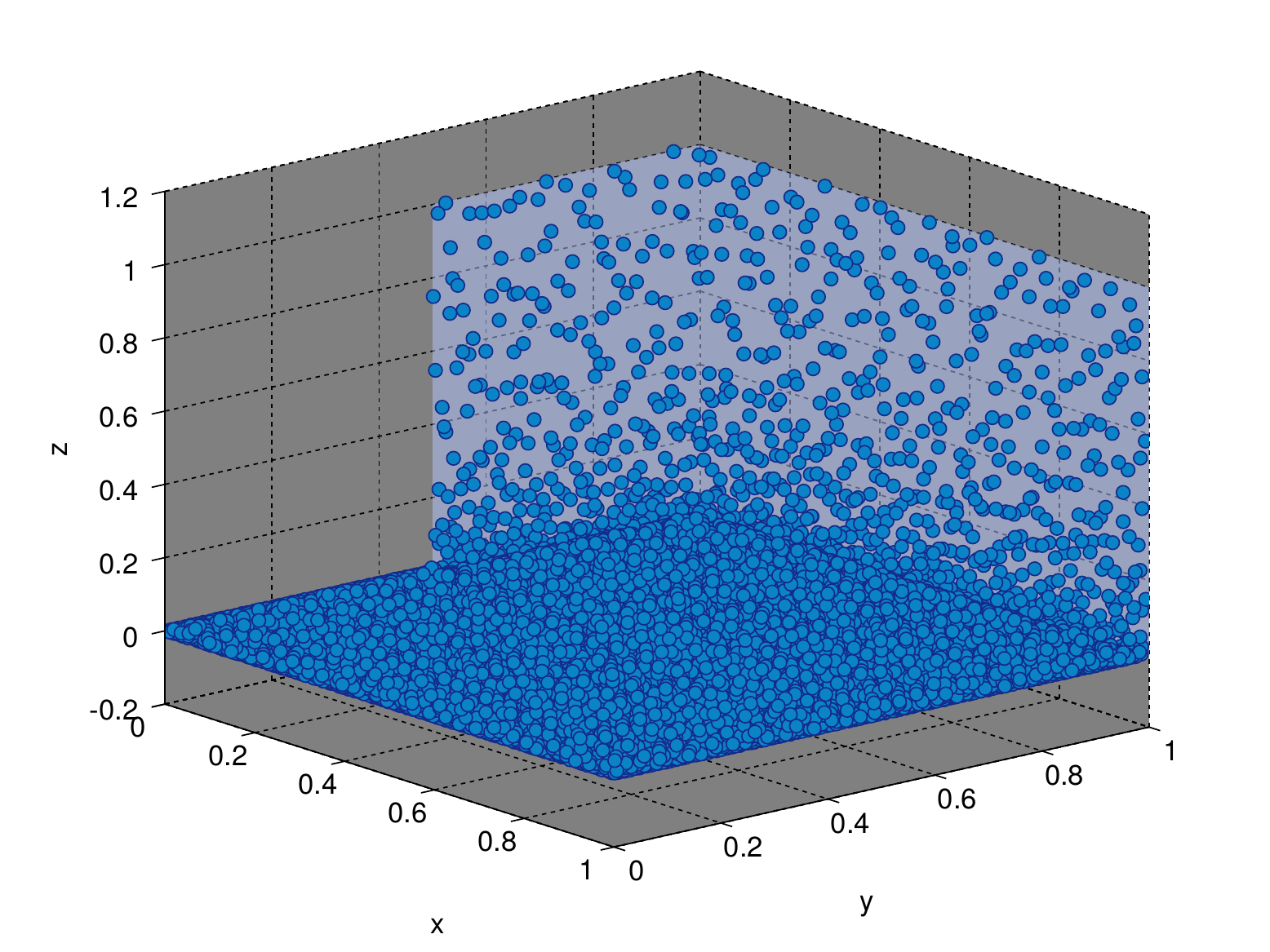}}\hspace{0.05\textwidth}
      \subfloat[Control, $c = 10^{-8}$]{\includegraphics[width=0.45\textwidth]{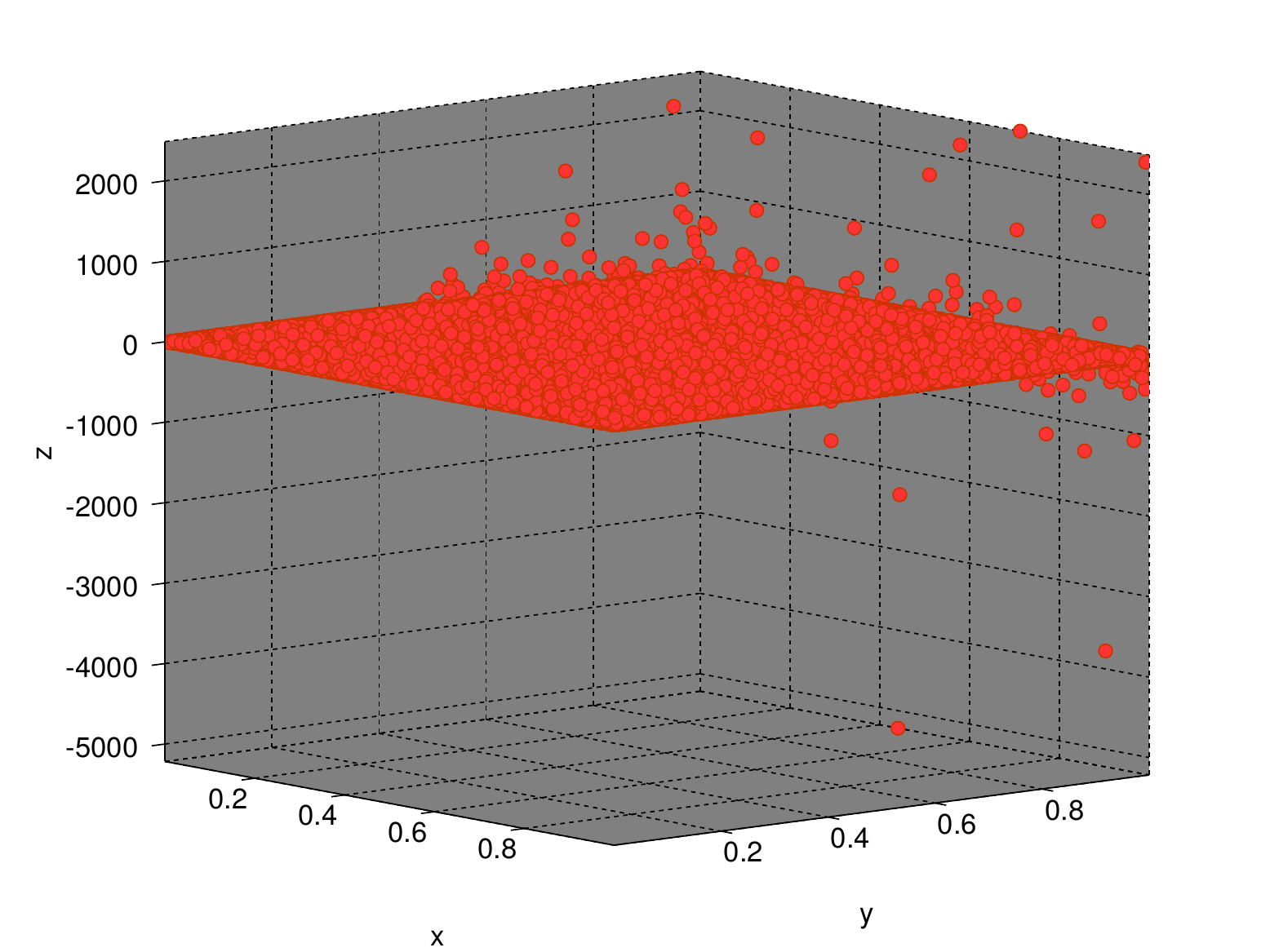}}
   \end{center}   
   \caption{Solution for problem 2, $n = 50000$, $\beta = 10^{-10}$}\label{figura:SolP6}
\end{figure}

	From the results reported in table \ref{tabla:ResP6} it can be seen that for any method it is possible to find $c$ in such a way that the minimum value of $\Vert y-\y\Vert_{L_{2}(\Omega)}$ is of the same order in magnitude. The difference is in $\Vert u\Vert_{L_{2}(\Omega)}$ since the norm obtained for the local scheme case is much smaller than for AC, which seems to affect in the same way the value of the cost functional.
	
	Figure \ref{figura:SolP6} only shows the solution for LAM-DQ using $n=50000$, since AC takes to much time, in fact, although we did not complete the experiment, we estimate that it will take around two days to obtain the results. The high number of total nodes were used to show the capabilities of LAM-DQ to handle big problems and to show in detail the solution obtained for the state for this problem. We can see how the solution is very close to $0$ in all the domain except very near of the boundary layer. 
	
\subsection{Problem 3}

	The last test problem is also a convection-diffusion control problem for which there is no exact solution, given by
\begin{align*}		
	\begin{array}{rrl}
		(-\epsilon\lap+\omega\cdot\nabla)y=u, & \beta(-\epsilon\lap-\omega\cdot\nabla)u= y-\y & \mbox{ in }\Omega         \\
		y=g,                                  & u=0                                           & \mbox{ on }\partial\Omega \\
		                                      &
	\end{array}\\
	\begin{array}{rll}
		& \y = \begin{cases}
			(2x_{1}-1)^{2}(2x_{2}-1)^{2} & \mbox{in }\left[0,\tfrac{1}{2}\right]^{2}\cap\Omega \\ 
			0                            & \mbox{elsewhere}
		  \end{cases}                                                                                     & \\
		& g = \begin{cases}
			(2x_{1}-1)^{2}(2x_{2}-1)^{2} & \mbox{in }\left[0,\tfrac{1}{2}\right]^{2}\cap\partial\Omega \\ 
			0                            & \mbox{elsewhere}
		  \end{cases}                                                                                     & \\
		&  \omega=(\cos{\theta},\sin{\theta}),\mbox{ with } \theta=2.4                                    & \\
		& \epsilon=\displaystyle\frac{1}{200}                                                             &
	\end{array}
\end{align*}

	Table \ref{tabla:ResP3} contains the same values as the previous example: $\Vert y-\y\Vert_{L_{2}(\Omega)}$ for the state $y$ and $\Vert u\Vert_{L_{2}(\Omega)}$ for control $u$.
{\scriptsize
	\begin{table}[htb]
		\begin{center}
			\begin{tabular}{c c c c c c c}\toprule[1.3pt]
				            & \multicolumn{3}{c}{\bf LAM-DQ}                      & \multicolumn{3}{c}{AC}\\ \cmidrule(lr){2-4} \cmidrule(lr){5-7}
				$\beta$     & $10^{-2}$       & $10^{-6}$       & $10^{-10}$      & $10^{-2}$       & $10^{-6}$       & $10^{-10}$      \\
				$c$         & \NSc{4.00}{-04} & \NSc{7.00}{-03} & \NSc{7.00}{-03} & \NSc{4.00}{-05} & \NSc{6.00}{-04} & \NSc{6.00}{-04} \\
				$y$         & \NSc{3.57}{-01} & \NSc{1.49}{-04} & \NSc{1.49}{-08} & \NSc{1.47}{-01} & \NSc{1.65}{-04} & \NSc{1.66}{-08} \\
				$u$         & \NSc{5.62}{00}  & \NSc{3.00}{00}  & \NSc{3.00}{00}  & \NSc{1.35}{00}  & \NSc{3.20}{00}  & \NSc{3.20}{00}  \\
				Cost        & \NSc{2.22}{-01} & \NSc{4.51}{-06} & \NSc{4.50}{-10} & \NSc{2.00}{-02} & \NSc{5.12}{-06} & \NSc{5.11}{-10} \\
				$\kappa$    & \NSc{6.52}{05}  & \NSc{1.25}{09}  & \NSc{6.52}{08}  & \NSc{2.63}{05}  & \NSc{2.41}{06}  & \NSc{2.41}{06}  \\
				$\kappa(S)$ & \NSc{5.25}{02}  & \NSc{2.27}{00}  & \NSc{1.00}{00}  &                 &                 &                 \\  
				Time        & \TimeFormat{00}{10} & \TimeFormat{00}{10} & \TimeFormat{00}{10} & \TimeFormat{01}{52} & \TimeFormat{01}{52} & \TimeFormat{01}{52} \\ \bottomrule[1.3pt]
			\end{tabular}\caption{Results from problem 3. For LAM-DQ $n^{(k)} = 50$. In both cases $n = 622$}\label{tabla:ResP3}
		\end{center}
	\end{table}
}
   
\begin{figure}[htb!]
   \begin{center}
      \subfloat{LAM-DQ, $n^{(k)} = 50$}\hspace{0.32\textwidth}
      \subfloat{AC\hspace{1.2cm}}\\
      \addtocounter{subfigure}{-2}
      \subfloat[State, $c = 10^{-6}$]{\includegraphics[width=0.41\textwidth]{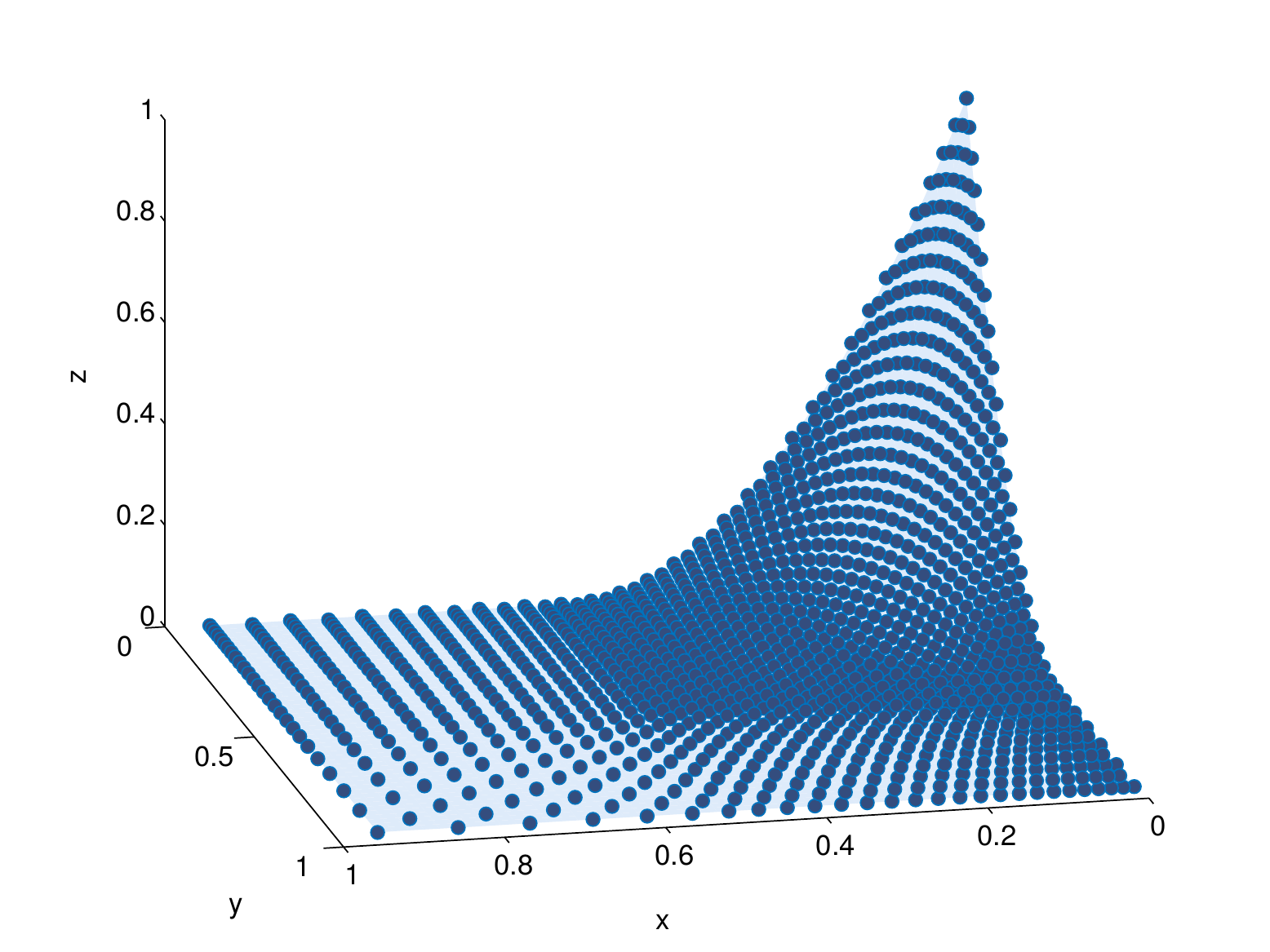}}\hspace{0.05\textwidth}
      \subfloat[State, $c = 10^{-7}$]{\includegraphics[width=0.41\textwidth]{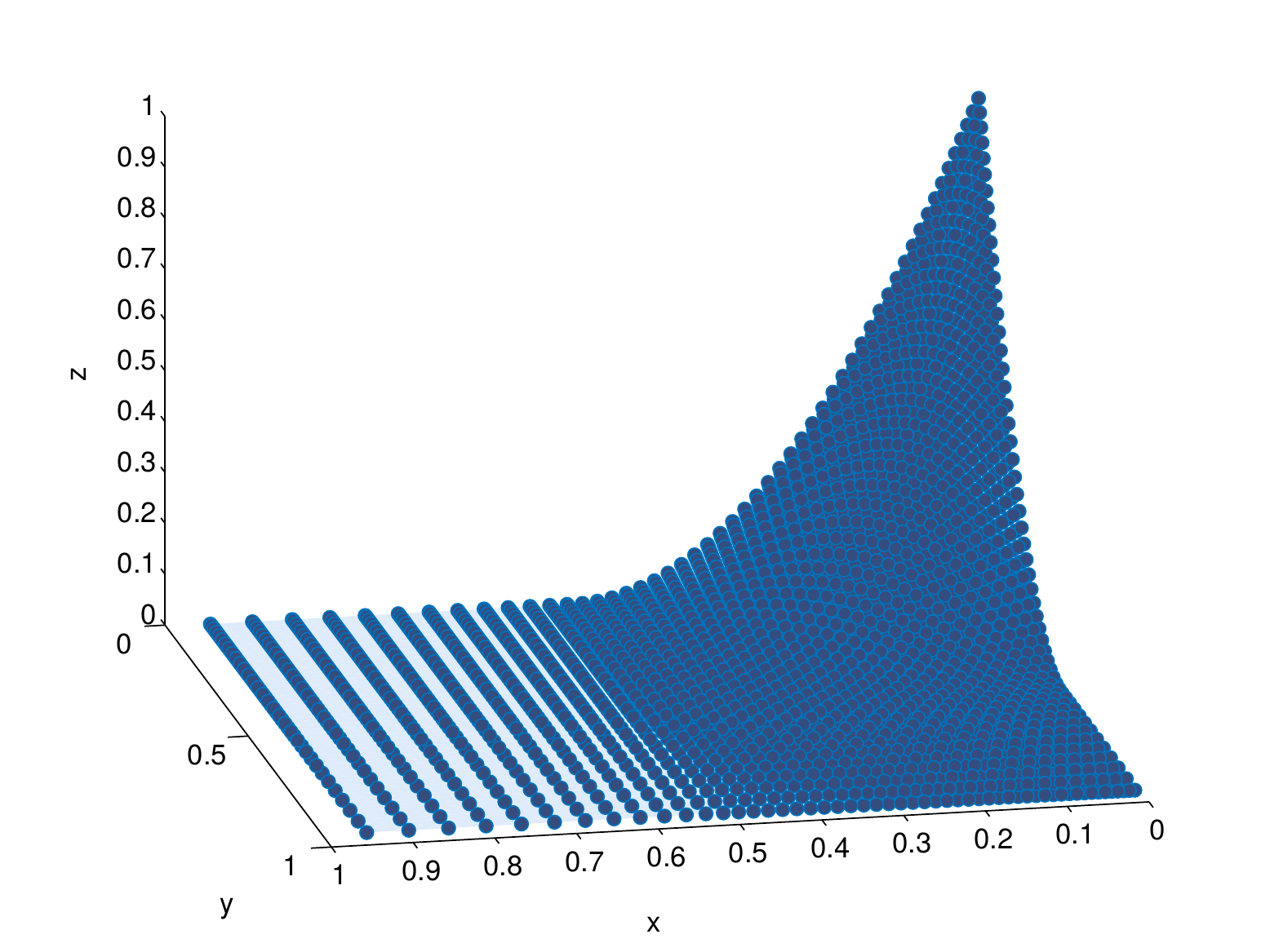}}\\
      \subfloat[Control, $c = 10^{-6}$]{\includegraphics[width=0.41\textwidth]{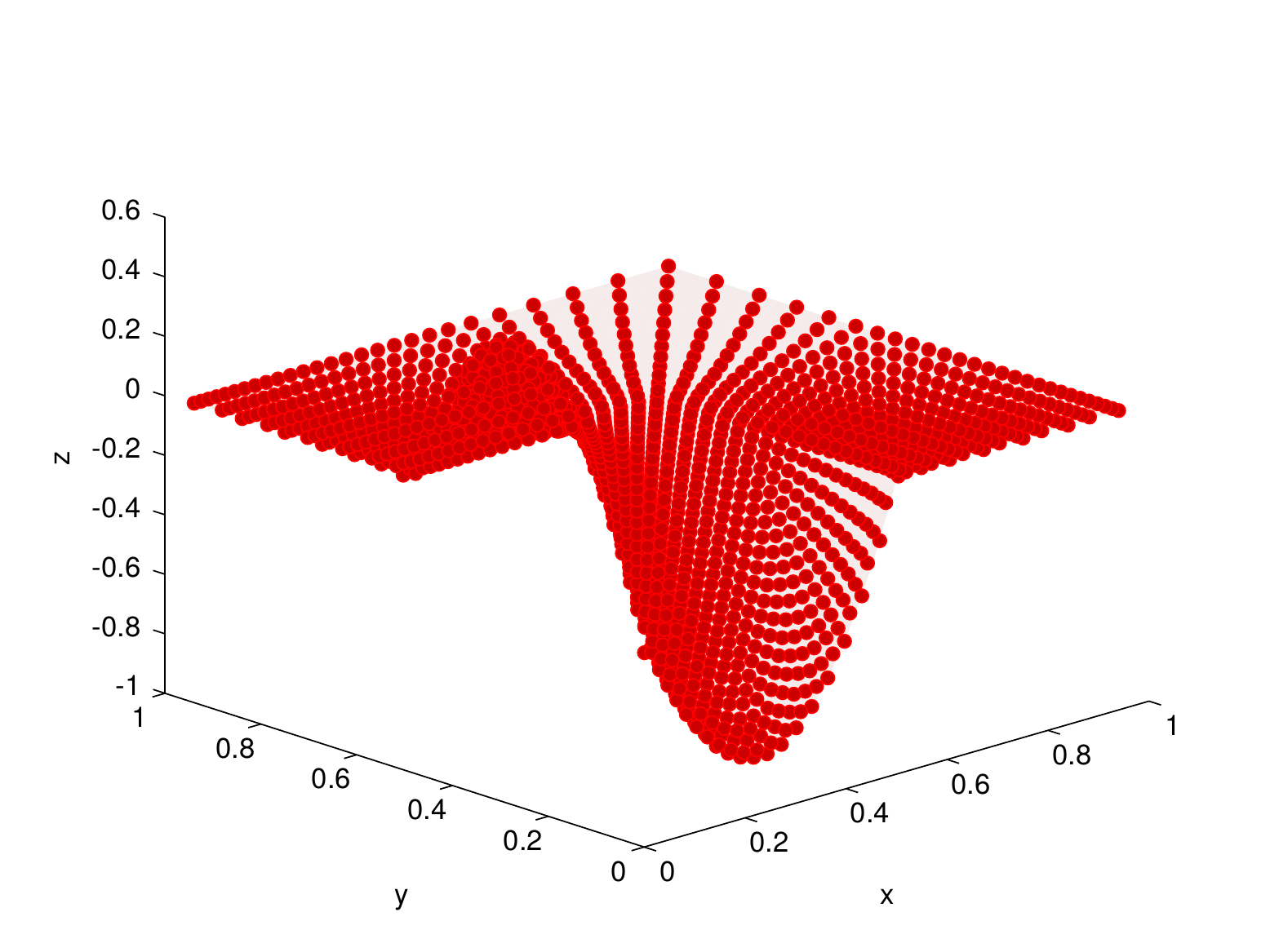}}\hspace{0.05\textwidth}
      \subfloat[Control, $c = 10^{-7}$]{\includegraphics[width=0.41\textwidth]{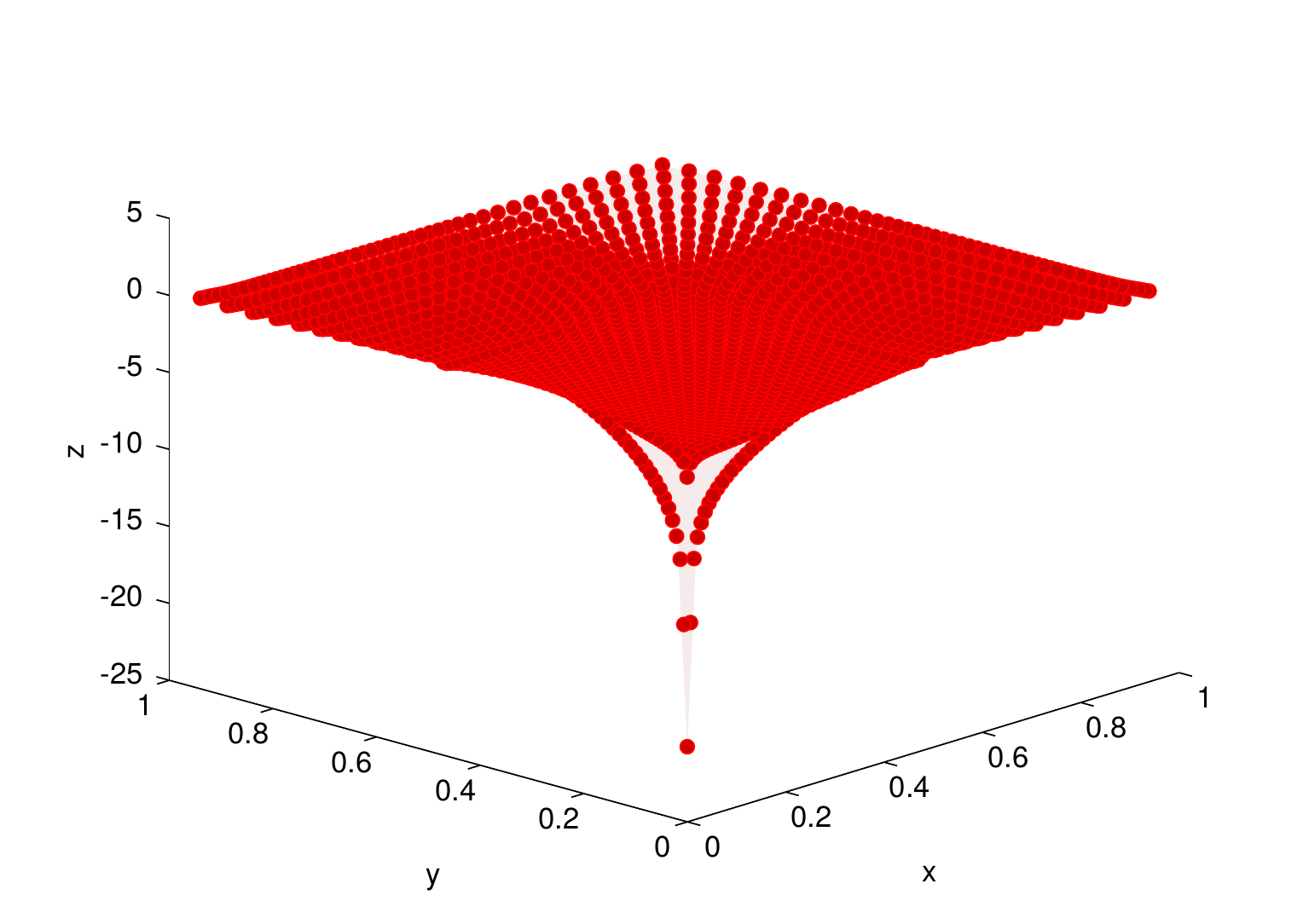}}
   \end{center}   
   \caption{Solution for problem 3, $n = 3021$, $\beta = 10^{-10}$}\label{figura:SolP3}
\end{figure}

	The results reported in the table \ref{tabla:ResP3} show again, that as in the previous example, for any method it is possible to find $c$ in such a way that the minimum value of $\Vert y-\y\Vert_{L_{2}(\Omega)}$ is of the same order in magnitude. Here, for the number of total nodes considered, there is no difference in the magnitude of $\Vert u\Vert_{L_{2}(\Omega)}$, and therefore also for the value of the cost functional.
	
	However, for the particular case for $\beta = 10^{-10}$ shown in figure \ref{figura:SolP3} with $n = 3021$, we have values for $y$ of the same magnitude but the control norm is lower for LAM-DQ. In addition, the control calculated through LAM-DQ visually resembles the results calculated by finite element method in \cite{Rees20101} and \cite{Rees20102}, which shows the consistency of the LAM-DQ solutions with respect to the finite element method.
	
\begin{figure}[t]
   \begin{center}
      \subfloat[Table]{{\scriptsize
			\begin{tabular}{c c c c}\toprule[1.3pt]
				\textbf{No. Nodos} & \textbf{AC}              & \textbf{LAM-DQ}     & \textbf{LAM-DQ Precond} \\ \cmidrule[1.1pt](lr){1-4}
				500                & \TimeFormat{00}{58}      & \TimeFormat{00}{07} & \TimeFormat{00}{15}     \\
				1000               & \TimeFormat{07}{34}      & \TimeFormat{00}{18} & \TimeFormat{00}{34}     \\
				1500               & \TimeFormat{24}{59}      & \TimeFormat{00}{36} & \TimeFormat{01}{01}     \\
				2000               & \TimeFormat{59}{09}      & \TimeFormat{01}{05} & \TimeFormat{01}{38}     \\
				2500               & \TimeFormat{114}{44}     & \TimeFormat{01}{49} & \TimeFormat{02}{29}     \\
				3000               & \TimeFormat{197}{50}     & \TimeFormat{02}{51} & \TimeFormat{03}{40}     \\
				3500               & \TimeFormat{316}{33}     & \TimeFormat{04}{15} & \TimeFormat{05}{19}     \\
  				4000               & \TimeFormat{467}{52}     & \TimeFormat{06}{06} & \TimeFormat{07}{12}     \\ \bottomrule[1.3pt]
			\end{tabular}}\label{tabla:ComTiempo}}\\
      \subfloat[Graph]{\includegraphics[width=0.55\textwidth]{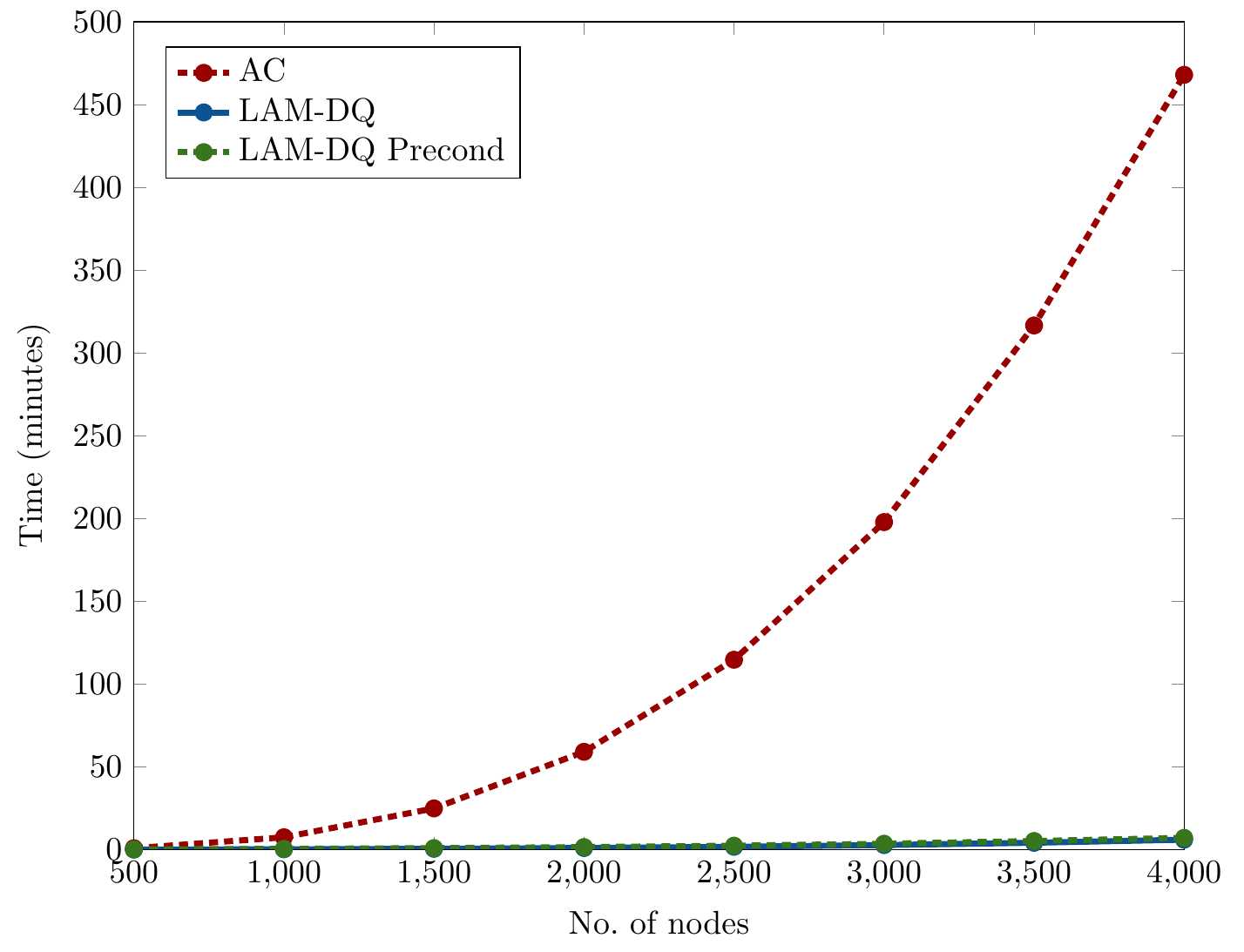}\label{figura:ComTiempo}}
   \end{center}
   \caption{Calculation time employed by the methods. For LAM-DQ, $n^{(k)}=50$.}
\end{figure}

	Finally, we compare the computing time for both methods. The tests were carried out using our own routines programmed in C++ on a machine with an Intel Core i5 M540 processor (2.53GHz). The execution time of the algorithms seems only to be dependent on the total number of nodes, that is, no matter which value of $c$ and $\beta$ are taken or if it is a convection-diffusion or Poisson control. Table \ref{tabla:ComTiempo} shows the different calculation times by varying the total number of nodes, showing that for all cases LAM-DQ has a shorter execution time in all cases. Figure \ref{figura:ComTiempo} shows in a more clear way the difference in the computing time between both methods, showing the advantage of LHI-DQ to solve massive problems.
   
\section{Conclusions}

	In this article, we solve control distributed problems for convection-diffusion linear PDEs problems by global and local radial basis functions methods. Inspired by the local Hermite interpolation method proposed by \cite{Stevens2010}, we introduced two new techniques, LAM-DQ and LAM-LAM.

A saddle point problem is obtained if we discretize the primal and adjoint equations by using radial basis function methods. We proposed a solution to this problem by discretizing instead, a well-posed biharmonic problem for the state variable and then obtaining the control by a second decoupled equation. 

An important contribution of this paper is that these local methods, in comparison to global colocation techniques, can attain similar precision errors for the same number of nodes, but with a considerable reduction of the computing, CPU, time.

While the condition number of the sparse global matrices in all our experiments, remains within an acceptable value, below the machine precision, the maximum condition number of the local matrices can grow up to the point where they are numerically singular as the fill distance tends to zero.

In this article, we deal with this problem by using quad precision and by proposing a simple but effective preconditioner.  By doing this, we manage to solve problems having 50000 nodes and reduce the condition number of the local matrices up to 10 orders of magnitude. The ill-conditioning of the Gram local and global matrices is currently an active research area in the field of radial basis function theory.

Although many proposals have appeared in the literature to deal with this general problem, we believe that the methods and the analysis proposed in this article present a significant contribution which shows the way to solve large distributed control problems.

\section*{Acknowledgements}
Funding: This work was supported by the National Autonomous University of M\'exico [grant: PAPIIT,  IN102116] 
\bibliography{Manuscript_arXiv}

\end{document}